\documentclass[11pt, final]{amsart}
\usepackage{amssymb}

\theoremstyle{plain}
\newtheorem{theorem}{Theorem}

\theoremstyle{remark}

\def\C1{C_1^+(H)}
\def\la{\langle}

\def\ra{\rangle}
\def\tr{\operatorname{tr}}

\begin{document}
\title[]{Fidelity preserving maps on density operators}
\author{LAJOS MOLN\'AR}
\address{Institute of Mathematics and Informatics\\
         University of Debrecen\\
         4010 Debrecen, P.O.Box 12, Hungary}
\email{molnarl@math.klte.hu}
\thanks{  This research was supported by the
          Hungarian National Foundation for Scientific Research
          (OTKA), Grant No. T030082, T031995, and by
          the Ministry of Education, Hungary, Reg.
          No. FKFP 0349/2000}
\date{\today}
\begin{abstract}
We prove that any bijective fidelity preserving transformation on the set of all
density operators on a Hilbert space is implemented by an either
unitary or antiunitary operator on the underlying Hilbert space.
\end{abstract}
\maketitle

Let $H$ be a Hilbert space. The set of all density operators on
$H$, that is, the set of all positive self-adjoint operators on
$H$ with finite trace is denoted by $\C1$. (We note that one may
prefer normalized density operators; see the first remark at the
end of the paper.)

According to Uhlmann \cite{Uhlmann,Uhlmann2}, for any $A,B \in
\C1$ we define the fidelity of $A$ and $B$ by
\[
F(A,B)=\tr (A^{1/2}B A^{1/2})^{1/2}.
\]
This is in fact the square-root of the transition probability
introduced by Uhlmann in \cite{Uhlmann3} for density operators
which later Jozsa called fidelity and showed its use in quantum
information theory \cite{Jozsa}. The reason that Uhlmann defined
the fidelity in the way above is that after taking square-root the
function $F$ behaves significantly better.

Clearly, the fidelity is in intimate connection with the
transition probability between pure states. Wigner's theorem
describing the form of all bijective transformations on the set of
all pure states which preserve the transition probability plays
fundamental role in the theory of quantum systems. By analogy, it
seems to be of some importance to describe all bijective
transformations on the density operators which preserve the
fidelity. This is exactly what we are performing in the present
paper. We show that the fidelity preserving transformations on
$\C1$ are implemented by an either unitary or antiunitary operator
of the underlying Hilbert space.

The main result of the paper reads as follows.

\begin{theorem}\label{T:fidel1}
Let $\phi :\C1 \to \C1$ be a bijective transformation with the
property that
\begin{equation}\label{E;fidel1}
F(\phi(A),\phi(B))=F(A,B) \qquad (A,B \in \C1).
\end{equation}
Then there is an either unitary or antiunitary operator $U:H\to H$
such that
\[
\phi(A)=UAU^* \qquad (A\in \C1).
\]
\end{theorem}

\begin{proof}
The main point of the proof is to reduce the problem to Wigner's
classical result. In order to do so, we first prove that $\phi$
preserves the order $\leq$ (which comes from the usual order
between bounded self-adjoint operators on $H$) on $\C1$. If
$A,B\in \C1$, $A\leq B$, then for any $C\in \C1$ we have
\[
C^{1/2}AC^{1/2}\leq C^{1/2}BC^{1/2}.
\]
Since the square-root function is operator monotone we have
\[
(C^{1/2}AC^{1/2})^{1/2}\leq (C^{1/2}BC^{1/2})^{1/2}.
\]
Taking trace we obtain
\[
F(A,C)\leq F(B,C)
\]
which implies
\[
F(\phi(A),\phi(C))\leq F(\phi(B),\phi(C))
\]
for every $C\in \C1$. Let $\phi(C)$ run through the set of all
rank-one projections. If $P$ is the rank-one projection projecting
onto the subspace generated by the unit vector $x\in H$, then we
have
\[
\la \phi(A)x,x\ra^{1/2}= F(\phi(A),P)\leq F(\phi(B),P)=\la
\phi(B)x,x\ra^{1/2}.
\]
Since this holds for every unit vector $x\in H$, we obtain
$\phi(A)\leq \phi(B)$. Since $\phi^{-1}$ has the same properties
as $\phi$, it follows that $\phi$ preserves the order in both
directions.

We next show that $\phi$ preserves the rank-one operators. In
fact, one can easily see that an element $A\in \C1$ is of rank one
if and only if the set $\{ T\in \C1 \, :\, T\leq A\}$ is infinite
and total in the sense that any two elements in it are comparable
with respect to the order $\leq $. By the order preserving
property of $\phi$ it now follows that $\phi$ preserves the
rank-one elements of $\C1$ in both directions.

Clearly, a rank-one operator $A\in \C1$ is a rank-one projection
if and only if its trace is 1, that is, if $F(A,A)=1$. It follows
that $\phi$ preserves the rank-one projections.

It needs elementary computation to show that for any rank-one
projections $P,Q$ we have
\[
F(P,Q)=(\tr PQ)^{1/2}.
\]
Hence, we have proved that if we restrict $\phi$ onto the set of
all rank-one projections, we have a bijective function on this set
which satisfies
\[
\tr \phi(P)\phi(Q)=\tr PQ
\]
for every $P,Q$. Now, we can apply Wigner's theorem and get that
there exists an either unitary or antiunitary operator $U:H\to H$
such that
\[
\phi(P)=UPU^*
\]
holds for every rank-one projection $P$. Replacing $\phi$ by the
transformation
\[
A\mapsto U^*\phi(A)U
\]
if necessary, we can obviously assume that our original
transformation $\phi$ satisfies $\phi(P)=P$ for every rank-one
projection $P$. It remains to show that $\phi(A)=A$ holds for
every density operator $A\in \C1$. If $x\in H$ is a unit vector
and $P$ is the corresponding rank-one projection, then we compute
\[
\begin{gathered}
\la \phi(A)x,x\ra^{1/2}=\tr (P\phi(A)P)^{1/2}=\\
F(\phi(A),P)= F(\phi(A), \phi(P))=F(A,P)=\la Ax,x\ra^{1/2}.
\end{gathered}
\]
Since this holds for every unit vector $x\in H$, we conclude that
$\phi(A)=A$ $(A\in \C1)$. This completes the proof.
\end{proof}

If the underlying Hilbert space is finite dimensional, then we can
get rid of the assumption on bijectivity and obtain our second
result which follows.

\begin{theorem}\label{T:fidel2}
Let $H$ be a finite dimensional Hilbert space and let $\phi:\C1\to
\C1$ be a transformation such that
\begin{equation}
F(\phi(A),\phi(B))=F(A,B) \qquad (A,B \in \C1).
\end{equation}
Then there is an either unitary or antiunitary operator $U:H\to H$
such that
\[
\phi(A)=UAU^* \qquad (A\in \C1).
\]
\end{theorem}

\begin{proof}
If $A,B$ are self-adjoint operators, then we say that $A,B$ are
mutually orthogonal if $AB=0$. Clearly, $A,B$ are mutually
orthogonal if and only if they have mutually orthogonal ranges.

Let us assume that the dimension of $H$ is $d\geq 2$. It is easy
to see that one can characterize the positive rank-one operators
in the following way: a positive operator $A$  is of rank one if
and only if $A\neq 0$ and there exists a system $A_1, \ldots,
A_{d-1}$ of nonzero positive operators such that the elements of
$A, A_1, \ldots, A_{d-1}$ are mutually orthogonal.

It is clear that $\phi$ preserves the nonzero operators. Indeed,
this follows form the equality $F(A,A)=\tr A$. Let $A,B$ be
positive operators with $AB=0$. Then $A,B$ are commuting and by
the properties of the positive square-root of positive operators,
we have the same for $A,B^{1/2}$. Therefore, we infer
\[
B^{1/2}AB^{1/2}=AB=0.
\]
So, we have $\tr (\phi(B)^{1/2}\phi(A)\phi(B)^{1/2})^{1/2}=0$. But
this implies that
\[
(\phi(B)^{1/2}\phi(A)\phi(B)^{1/2})^{1/2}=0.
\]
Hence, we have
\[
\phi(B)^{1/2}\phi(A)\phi(B)^{1/2}=0
\]
from which we get
\[
(\phi(A)^{1/2}\phi(B))^*(\phi(A)^{1/2}\phi(B))=
\phi(B)\phi(A)\phi(B)=0.
\]
Consequently, we have $\phi(A)^{1/2}\phi(B)=0$ which implies
$\phi(A)\phi(B)=0$. This shows that $\phi$ preserves the
orthogonality in one direction. By the characterization of
rank-one operators given in the beginning of the proof, we infer
that $\phi$ sends rank-one operators to rank-one operators. Now,
similarly to the corresponding part of the proof of our previous
theorem one can check that $\phi$ sends rank-one projections to
rank-one projections. Just in that proof one can readily verify
that
\[
\tr \phi(P)\phi(Q)=\tr PQ
\]
holds for every rank-one projection $P,Q$. Now, by the
nonsurjective version of Wigner's theorem \cite{Bargmann} (also
see \cite[Theorem 3]{Molnar}), we have an isometry or antiisometry
$U:H\to H$ such that
\[
\phi(P)=UPU^*
\]
holds for every rank-one projection $P$. Since $H$ is finite
dimensional, $U$ is in fact a unitary or antiunitary operator. The
proof can now be completed very similarly to the proof of our
first theorem.
\end{proof}

We conclude the paper with some remarks.

(1) In the introduction we have mentioned that one may prefer to
restrict the investigation to normalized density operators, that
is, to positive self-adjoint operators with trace 1. Although
following Uhlmann, in our treatment we did not insist on
normalization, we point out to the fact that one can get the same
result also in that case. The only thing that should be done is
the following. If $\phi$ is a bijective transformation on the set
of all normalized density operators, then we define $\tilde \phi
:\C1\to \C1$ in the following way: $\tilde \phi(0)=0$ and for any
$0\neq A\in \C1$ we set
\[
\tilde \phi(A)=\tr A \phi\biggl(\frac{A}{\tr A}\biggr).
\]
It is apparent that $\tilde \phi:\C1 \to \C1$ is a bijective
transformation extending $\phi$ and it preserves the fidelity.
Now, Theorem~\ref{T:fidel1} applies.

(2) One can easily generalize our results to obtain the same
description of transformations on density operators preserving not
the ''full'' fidelity but a certain part of it. We mean the
quantity $F_m^+(A,B)$ denoting the sum of the $m$ largest
eigenvalues of the operator $(A^{1/2}BA^{1/2})^{1/2}$ $(A,B\in
\C1)$ \cite[Definition]{Uhlmann2}. Here $m$ is fixed and when we
speak about eigenvalues we always take into account the
multiplicities. Now, one can formulate the same assertions as in
Theorem~\ref{T:fidel1} and \ref{T:fidel2} with $F_m^+$ in the
place of $F$. As for the proofs, one can follow quite the same
argument. In fact, the only additional thing that should be
observed concerns the order preserving property. Namely, one
should verify that $A\leq B$ if and only if $F_m^+(A,C)\leq
F_m^+(B,C)$ holds for every density operator $C$. The sufficiency
is clear if $C$ runs through the set of all rank-one projections.
The necessity follows from Weyl's monotonicity theorem stating
that if $A\leq B$, then the the $k$th largest eigenvalue of $A$ is
less than or equal to the $k$th largest eigenvalue of $B$ (cf.
\cite[Lemma 1.1, p. 26]{Gohberg}).

It would certainly be of interest to obtain similar results
concerning the ''partial'' fidelities introduced by Uhlmann in
\cite{Uhlmann}.

(3) It is easy to see that just as in Wigner's classical theorem,
the implementing unitary or antiunitary operator $U$ is unique up
to a phase factor (a scalar of modulus one).

(4) Following the lines in the proof of our first result one can
easily see that there is no need to assume the injectivity of the
transformation $\phi$. We set this condition only for the sake of
''symmetry'' in the formulation.

\vskip 10pt {\bf Acknowledgements.} The author is very grateful to
Armin Uhlmann for fruitful discussions on the topic of the paper
and also for encouraging the publication of the obtained results.

\bibliographystyle{amsplain}

\end{document}